\documentclass[12pt]{amsart}
\usepackage{txfonts,fancyvrb}

\usepackage{algorithm}
\usepackage[noend]{algpseudocode}
\makeatletter
\def\BState{\State\hskip-\ALG@thistlm}
\makeatother

\usepackage{subfig}
\usepackage{amscd,amssymb,mathabx}
\usepackage[arrow,matrix,graph,frame,poly,arc,tips]{xy}
\usepackage{graphicx,calrsfs}
\usepackage{mathtools}
\usepackage[dvipsnames]{xcolor}
\usepackage{tikz}
\usetikzlibrary{decorations.markings}
\usetikzlibrary{shapes}
\usetikzlibrary{arrows}
\usetikzlibrary{backgrounds}
\usetikzlibrary{calc}
\usepackage{mdwlist}
\usepackage{multirow}
\usepackage{enumerate}
\usepackage{verbatim}
\usepackage{etoolbox}
\AtEndEnvironment{proof}{\setcounter{claim}{0}}

\newcommand\tto{\twoheadrightarrow}

\DeclarePairedDelimiter{\form}{\langle}{\rangle}

\DeclarePairedDelimiter{\floor}{\lfloor}{\rfloor}

\newcommand\ba{\begin{align*}}
\newcommand\ea{\end{align*}}
\newcommand\be{\begin{enumerate}}
\newcommand\ee{\end{enumerate}}
\newcommand\bp{\begin{proof}}
\newcommand\ep{\end{proof}}
\newcommand\bpp{\begin{prop}}
\newcommand\epp{\end{prop}}
\newcommand\bpb{\begin{prob}}
\newcommand\epb{\end{prob}}
\newcommand\bd{\begin{defn}}
\newcommand\ed{\end{defn}}
\newcommand\bh{\begin{hint}}
\newcommand\eh{\end{hint}}

\newcommand\bC{\mathbb{C}}

\newcommand\bN{\mathbb{N}}
\newcommand\N{\mathbb{N}}
\newcommand\bR{\mathbb{R}}
\newcommand\R{\mathbb{R}}
\newcommand\bQ{\mathbb{Q}}
\newcommand\Q{\mathbb{Q}}
\newcommand\bZ{\mathbb{Z}}

\newcommand\FF{\mathcal{F}}

\newcommand\SL{\operatorname{SL}}

\newcommand\PSL{\operatorname{PSL}}

\newcommand\sse{\subseteq}
\newcommand\co{\colon}

\renewcommand{\MR}[1]
{\href{http://www.ams.org/mathscinet-getitem?mr=#1}{MR#1}}

\def\thetitle{{Non-freeness of groups generated by two parabolic elements with small rational parameters}}
\def\theauthors{{Sang-hyun Kim and Thomas Koberda}}
\usepackage{hyperref}
\hypersetup{
   colorlinks=false,
   plainpages,
   urlcolor=black,
   linkcolor=black
   pdftitle=   \thetitle,
   pdfauthor=  {\theauthors}
}

\theoremstyle{plain}
\newtheorem{thm}{Theorem}[section]
\newtheorem{lem}[thm]{Lemma}
\newtheorem{cor}[thm]{Corollary}
\newtheorem{prop}[thm]{Proposition}
\newtheorem{con}[thm]{Conjecture}
\newtheorem{que}[thm]{Question}
\newtheorem*{claim*}{Claim}

\newtheorem{claim}{Claim}

\newtheorem*{maincon}{\bf Main Conjecture}

\theoremstyle{remark}
\newtheorem{exmp}[thm]{Example}
\newtheorem{rem}[thm]{Remark}

\theoremstyle{definition}
\newtheorem{defn}[thm]{Definition}
\newtheorem{prob}{Problem}[section]

\begin{document}
\title\thetitle
\keywords{Fuchsian groups, Kleinian groups, Schottky groups, M\"obius groups}
\subjclass[2010]{Primary: 30F35, 30F40; Secondary: 20E05,  11J70}



\author[S. Kim]{Sang-hyun Kim}
\address{School of Mathematics, Korea Institute for Advanced Study (KIAS), Seoul, 02455, Korea}
\email{skim.math@gmail.com}
\urladdr{http://cayley.kr}

\author[T. Koberda]{Thomas Koberda}
\address{Department of Mathematics, University of Virginia, Charlottesville, VA 22904-4137, USA}
\email{thomas.koberda@gmail.com}
\urladdr{http://faculty.virginia.edu/Koberda}

\begin{abstract}
Let $q\in\bC$, let \[a=\begin{pmatrix} 1&0\\1&1\end{pmatrix},\quad
b_q=\begin{pmatrix} 1&q\\0&1\end{pmatrix},\] and let $G_q<\mathrm{SL}_2(\bC)$ be the group generated by $a$ and $b_q$. In this paper, we study the problem of determining when the group $G_q$ is not free for $|q|<4$ rational. We give a robust computational criterion which allows us to prove that if $q=s/r$ for $|s|\leq 27$ then $G_q$ is non-free, with the possible exception of $s=24$. In this latter case, we prove that the set of denominators $r\in\bN$ for which $G_{24/r}$ is non-free has natural density $1$. For a general numerator $s>27$, 
we prove that the lower density of denominators $r\in \bN$ for which $G_{s/r}$ is non-free has a lower bound 
\[
1- \left(1-\frac{11}{s}\right)
\prod_{n=1}^\infty \left(1-\frac{4}{s^{2^n-1}}\right).
\]
Finally, we show that for a fixed $s$, there are arbitrarily long sequences of consecutive denominators $r$ such that $G_{s/r}$ is non-free. The proofs of some of the results are computer assisted, and Mathematica code has been provided together with suitable documentation.
\end{abstract}

\maketitle
\setcounter{tocdepth}{1}
\section{Introduction}
For each $q\in\bC$, let us write
\[a=\begin{pmatrix} 1&0\\1&1\end{pmatrix},\quad
b_q=\begin{pmatrix} 1&q\\0&1\end{pmatrix},\] and write $G_q$ for the subgroup of $\mathrm{SL}_2(\bC)$ generated by $a$ and $b_q$.

The group $G_q$ is not cyclic unless $q=0$. 
It is proved by Sanov~\cite{Sanov47} (for $|q|=4$) and Brenner~\cite{Brenner1955} (for $|q|>4$) that the group $G_q$ is free for all $q\in\bR\setminus(-4,4)$; more strongly, the group $G_q$ is discrete and free for all $q$ in the \emph{Riley slice} of the complex plane~\cite{KS1994}.

In this paper, we study the following conjecture:
\begin{maincon}\label{con:main}
For each nonzero rational number $q=s/r$ in $(-4,4)$, the group
\[
G_q:=\form{a,b_q}\le\SL_2(\bC)\]
is not free.
\end{maincon}

Lyndon and Ullman asked this conjecture (as a question) in~\cite{LyndonUllman69}.
This problem has a long history, and the reader is directed to~\cite{Gilman2008} and to Section~\ref{sec:notes} below for the state of the art prior to this writing. 

Slightly different normalizations have also been considered in the literature. 
We may define 
 \[H_q=\form*{
 \begin{pmatrix} 1&0\\q&1\end{pmatrix},\quad
\begin{pmatrix} 1&q\\0&1\end{pmatrix}}.\] 
The corresponding question for $H_q$ is attributed to Merzlyakov in the Kourovka Notebook~\cite[Problem 15.83]{Kourovka2018}. 
It is noted in~\cite{CJR58} that  $H_q\cong G_{q^2}$.
In some other papers such as~\cite{CJR58,Gilman2008}, the group $G_{2q}$ is considered.

\begin{rem} Under the hypothesis that $q$ is rational and belongs to $(-4,4)$, the group $G_q$ is discrete only if $|q|\in\{0,1,2,3\}$; see~\cite{Knapp1969}.\end{rem}

\subsection{Main results}
As mentioned above, $G_q$ is free whenever $q\in\bR\setminus(-4,4)$. It is easy to see that $G_q$ is free if $q$ is transcendental.
However, being algebraic is not sufficient to guarantee non-freeness.
As noted in~\cite{CGS2017},
Galois conjugation yields an isomorphism 
\[
G_{4-\sqrt{2}}\cong G_{4+\sqrt{2}},\]
the latter of which is indeed free by the result of Sanov and Brenner.

\bd\label{d:relation number}We will say $q\in\bC$ is a \emph{relation number} if $G_q$ is not a rank--two free group.\ed

A good summary of known results about \emph{rational} relation numbers can be found in Theorem 7.7 of ~\cite{Gilman2008}.
Before stating the results of this paper, we introduce some terminology.
Let $F=\langle x,y\rangle$ be a free group of rank two.
A complex number $q$ is called an \emph{$\ell$--step relation number} 
if there exists a nontrivial word of the form
\[
w=y^{m_1}x^{m_2}\cdots y^{m_{2k+1}}\in F\]
for some $k\in[0, \ell]$ and $m_i\in\bZ\setminus\{0\}$ such that $w(a,b_q)$ is a lower--triangular matrix in $\SL_2(\bC)$.

It turns out then every relation number is an $\ell$--step relation number for some $\ell\ge0$, and vice versa (Lemma~\ref{l:equiv-odd-word}).
Actually, if $q$ is an $\ell$--step relation number, then there exists a word $v=v(x,y)\in F$ of syllable length at most $8(\ell+1)$ such that $v(a,b_q)=1$; see Remark~\ref{r:equiv-odd-word}.

Let $X\sse\bZ$ be a subset. 
The \emph{(right) upper density} of $X$ is given by \[\overline{d}(X)=\limsup_N \frac{X\cap [1,N]}{N}.\] The
\emph{(right) lower density} of $X$ is similarly given by  \[\underline{d}(X)=\liminf_N \frac{X\cap [1,N]}{N}.\] If these limits coincide, they are
called the \emph{(right) natural density} of $X$.
Note we allow $X$ to have negative integers. 
\begin{rem}
One may also consider a \emph{symmetric (lower or upper) density}, which is a limit (superior or inferior) of $(X\cap[-N,N])/(2N+1)$. 
For the integer sets concerned in this paper, all right densities will coincide with symmetric densities, whence we will simply refer to upper and lower densities when no confusion can arise. Note in particular that if $s/r$ is an $\ell$--step relation number then so is $s/(-r)$.\end{rem}

Our main results are towards resolving the Main Conjecture. Precisely, we prove the following:

\begin{thm}\label{thm:comput}
Let $s$ be a positive integer.
\be
\item
Suppose $s\leq 27$ and $s\ne 24$.
Then for all but finitely many nonzero integers $r$, the number $s/r$ is a $2$--step relation number.
Moreover, for all nonzero integer $r$ satisfying $s/r\in(-4,4)$,
the number $s/r$ is a relation number.
\item\label{p:s24}
If $s=24$, then $s/r$ is a $2$--step relation number
for all $r$ in some natural density--one subset of $\bN$.
\ee
\end{thm}

By our previous discussion, the above theorem resolves the Main Conjecture for all $r$ if $s\in[1,27]\setminus\{24\}$, and for almost all $r$ if $s=24$.
It even asserts that for a given $s\le27$ and for almost all $r\in\bN$, there exists a nontrivial word of syllable length at most $24$ in $G_q$ that becomes trivial.
We note that some parts of the proof are computer assisted, and we have provided code and documentation in the appendices below.

For a general $s\in\N$, we have the following result which finds a very large number of relation numbers with a given numerator:

\begin{thm}\label{thm:density}
Let $s$ be an integer greater than $27$. If we set
\[
A_s^{(2)}:=\{r\in\bZ\setminus\{0\}\mid s/r\text{ is a $2$--step relation number}\},\]
then we have
\[ \underline{d}\left(A_s^{(2)}\right)\ge 1- \left(1-\frac{11}{s}\right)\prod_{n=1}^\infty \left(1-\frac{4}{s^{2^n-1}}\right).\]
\end{thm}

It is natural to wonder if ${d}\left(A_s^{(2)}\right)=1$. Unfortunately, the sequence $\{s^{2^n-1}\}_{i\in\N}$ grows much too quickly, and generally
the infinite product in Theorem~\ref{thm:density} will converge to real number strictly less than $1$ (see Section~\ref{sec:gde} below). Of course, the choices of such a sequence can be modified, but it is not clear to the authors that the methods given here avail themselves to a suitable choice that
witnesses ${d}\left(A_s^{(2)}\right)=1$. 

\begin{que}\label{que:zero-dense}
For an integer $s>27$, is it true that ${d}\left(A_s^{(2)}\right)=1$?
\end{que}

We are able to prove one further result which strongly suggests that the answer to Question~\ref{que:zero-dense} is yes, without quite
establishing it definitively.


\begin{thm}\label{thm:progression}
(see Corollary~\ref{cor:progression})
Let $s,r,N\in\bN$.
Then there exists an $M=M(s,r,N)\in \bN$ such that 
\[
\frac{s}{r+i+sMj}\]
is a $2$--step relation number for all integers $0\le i< N$ and $j\ne0$.
\end{thm}

In particular, for a fixed $s$ there are arbitrarily long sequences of consecutive denominators which give rise to relation numbers of the form $s/r$. However, such sequences may possibly be spaced very sparsely within $\bN$.

\subsection{Notes and references}\label{sec:notes}
As noted above, the extent to which Sanov's result holds or fails for $q\in (-4,4)$ has a long history.
Some of the earliest examples of non-integral rational relation numbers of $q$ were found by Ree~\cite{Ree61}. On the other hand, many conditions for freeness of $G_q$ were found by
Chang--Jennings--Ree~\cite{CJR58}. 
Many more examples of relation numbers were found in~\cite{BMO75,Ignatov86,Ignatov90,LyndonUllman69,Bamberg2000}. 
Connections to diophantine problems, and especially
solutions to Pell's Equation, were studied in~\cite{Farbman95,TanTan96,Beardon2003}.
Discreteness of $G_q$ for a complex parameter $q\in\bC$ has been extensively studied; see~\cite{ASWY2007,Gilman2008} and the references therein. For related discreteness questions in $\PSL_2(\R)$, see~\cite{KKM2019}, for instance.

A dynamical interpretation of relation numbers was suggested first by Tan--Tan~\cite{TanTan96}, and these ideas have been developed
in~\cite{Bamberg2000,Romik2008,Slanina2015}.

One may compare the results of this paper to the results outlined in Theorem 7.7 of~\cite{Gilman2008}. We are primarily concerned with
groups of the form $G_q$ for $|q|<1$ rational, whereas the results there are given for groups of the form $H_q$ where $q$ may be
non-rational algebraic. One notes immediately from Theorem~\ref{thm:comput} that we have produced many new examples of rational
relation values of $q$, and in view of Theorems~\ref{thm:density} and~\ref{thm:progression}, many new infinite families of relation
values which do not fall under the purview of previously known results.

The freeness and non-freeness of the groups $G_q$ has applications to group--based cryptography and theoretical computer science. See for
instance~\cite{CGS2017}.

Finally, a remark about normalizations. We consider the groups $\{G_q\}_{q\in\Q}$ over the groups $\{H_q\}_{q\in\Q}$, in spite of the break
in symmetry, because the groups $\{G_q\}_{q\in\Q}$ encompass
a larger class of subgroups of $\SL_2(\Q)$ and hence give rise to an \emph{a priori} richer theory.

\subsection{General strategy and intuition}

Our approach to Conjecture~\ref{con:main} is essentially from first principles. If $q$ is a parameter for which $G_q$ is not free then
very elementary manipulations show that $q$ has to be a root of a polynomial with rational coefficients. The degrees of these polynomials
are related to the simplest nontrivial words in the free group which witness the fact that $G_q$ is not free, where here complexity is
measured in terms of the \emph{syllable length} of words.

For high degree polynomials, criteria for defining natural families of relation numbers are difficult to formulate in a way which is concise and
amenable to study, so that we restrict our attention to relatively simple polynomials. From there, we consider the following question:
what conditions on $r\in\bZ$ force $s/r$ to be a relation number for $s\in\bZ$ fixed?

The answers we propose have to do with the divisibility properties of $r$ modulo various multiples of $s$. This leads to many  technical
definitions (cf. $s$--good residue classes in Definition~\ref{d:sgood} below), and the main technical tools (see Lemmata
~\ref{l:quad0}, ~\ref{l:quad1}. and~\ref{l:quad2} below). These tools allow us to declare
all sufficiently large elements of certain residue classes modulo some multiple of $s$ to be relation numbers.

So, to show that $s/r$ is always a relation number for fixed $s$ and $r> s/4$, we begin showing as many residue classes as possible modulo
$sm$ consist of relation numbers, for some nonzero integer $m$. Then, take the remaining residue classes and consider their
residues modulo $sm'$ for some multiple $m'$ of $m$. Then, the technical tools allow us to conclude that many of these residue classes
modulo $sm'$ consist of relation numbers. Through this recursive procedure, more and more values of $r$ are shown to give relation numbers,
and the hope is that the procedure terminates after finitely many steps.

For $s\leq 27$ and $s\neq 24$, we can indeed show that the procedure terminates in finitely many steps, proving that all the relevant
rational parameters with those numerators are relation numbers. For $s=24$, we cannot show that the procedure terminates, but we have
enough control over the number of residue classes which are eliminated at each stage to conclude that the set of denominators for which
$24/r$ is not a relation number has natural density zero. We generalize these ideas to give lower bounds on the natural density of
relation number denominators for arbitrary numerators.

\section{Notation and terminology}
Recall we have separately defined a \emph{relation number} and an \emph{$\ell$--step relation number} in the introduction.
The number $q=0$ is the unique $0$--step relation number. 
The following lemma (due to Lyndon and Ullman) describes the relationship between the Main Conjecture and $\ell$--step relation numbers.


\begin{lem}[\cite{LyndonUllman69}]\label{l:equiv-odd-word}
A complex number $q$ is a relation number if and only if
it is an $\ell$--step relation number for some $\ell\ge0$.\end{lem}
\bp The forward direction is obvious from the fact that the identity matrix is lower--triangular.  For the converse, let $w=w(x,y)$ be
such that the matrix
\[w(a,b_q)\cdot a\cdot w(a,b_q)^{-1}\]
is lower triangular such that the diagonal entries are $1$.
It follows that the reduced word $[wxw^{-1},x]$
  becomes the identity in $\SL_2(\bC)$ after setting $x=a$ and $y=b_q$.\ep
\begin{rem}\label{r:equiv-odd-word}
The \emph{syllable length}  of  a nontrivial element $g\in F$ is the smallest integer $\ell\ge0$ such that  \[ g=w_1\cdots w_\ell\] for some $w_i\in \form{x}\cup \form{y}$.
The above proof shows that if $q$ is an $\ell$--step relation number then there exists a nontrivial word $v(x,y)=[wxw^{-1},x]$ of syllable length at most $8(\ell+1)$ such that $v(a,b_q)=1$.\end{rem}

From Lemma~\ref{l:equiv-odd-word},  we see that the Main Conjecture has the following diophantine-type formulation.
\begin{con}\label{con:main-hi}
Every rational number in  $(-4,4)$ is an $\ell$--step relation number for some $\ell\ge0$.
\end{con}

Let us describe a notation that will be used often throughout this paper. 
Let $q\in\bC$, and let $m_1,m_2,\ldots$ be a sequence of nonzero integers.
We define complex vectors $v_1,v_2,\ldots$ by setting $v_1=(1,0)$ and 
\[v_{i+1} = (1, 0)b_q^{m_1}a^{m_2}\cdots (b_q\text{ or }a)^{m_i}.\]
Note that $q$ is an $\ell$--step relation number if and only if one can find a sequence $\{m_i\}\sse\bZ\setminus\{0\}$ such that $v_{2k+2}\in\bC\times \{0\}$ for some $k\le\ell$.

As we are only interested in whether or not the second coordinate of $v_i$ becoming zero, we may regard $v_i$ as a point in the projective space $\bC P^1$. 
In particular, we will identify $(x,y)$ and $(nx,ny)$ for $x,y\in\bZ$ and $n\in\bZ\setminus\{0\}$. We will then use the notation
\begin{equation}\tag{$*$}\label{eq:orb}
v_1:=(1,0)\stackrel{m_1}{\tto}v_2\stackrel{m_2}{\to}v_3\stackrel{m_3}{\tto}\dots\stackrel{m_{2i}}{\to}v_{2i+1}\tto\cdots.\end{equation}
The nonzero exponents $m_1,m_2,\ldots$ will often be suppressed as well.
\begin{exmp}\label{exmp:int}
For $q=1$ or $q=2$, we have a sequence 
\[
(1,0)\tto (1,2)\to(-1,2)\tto(-1,0)=(1,0).\]
For $q=3$, we see
\[
(1,0)\stackrel{1}{\tto}(1,3)\stackrel{-1}{\to} (-2,3)\stackrel{1}{\tto}(-2,-3)\stackrel{-1}\to(1,-3)\stackrel{1}\tto(1,0).\]
It follows that all integers in the interval $[-3,3]$ are relation numbers.
\end{exmp}

The Main Conjecture can be reformulated in terms of generalized continued fractions. Suppose we have an orbit as above in (\ref{eq:orb}).
Write $Q=1/q$ and $v_i=(x_i,y_i)$. Assuming $x_i y_i\ne0$, we define
\[q_i:=Qy_i/x_i=y_i/(q x_i).\] 
Then we have that
\[
q_{i+1}=\begin{cases}
{Q(y_i+qm_i x_i)}/{x_i}=m_i+q_i&\text{ if }2\not\divides i,\\
{Qy_i}/({x_i+m_i y_i})=Q/(m_i+Q/q_i)&\text{ if }2\divides i.
\end{cases}
\]

On the other hand, it is obvious that $q$ is a relation number if 
$x_iy_i=0$ for some $i$, or if \[(x_i,y_i)=(x_j,y_j)\in \bC P^1\] for some $i<j$.
In summary, we have the following.
\begin{prop}\label{p:CF}
Let $Q\in\bC\setminus\{0\}$. Then $1/Q$ is a relation number if and only if 
there exists a finite sequence of non-zero integers
\[m_1,\ldots,m_\ell\]
such that the sequence 
\[
a_k:= m_k+\cfrac{Q}{m_{k-1}+\cfrac{Q}{\cdots+\cfrac{Q}{m_2 +\cfrac{Q}{m_1}}}}\]
either terminates with $a_\ell=0$ for some $\ell\ge2$, or 
satisfies $a_\ell=a_{\ell'}$ for some $\ell>\ell'\ge 2$.
\end{prop}
The Main Conjecture asserts that one has a sequence $\{m_i\}$ as above whenever $Q$ is a rational number satisfying $|Q|>1/4$.

\section{Families of rational relation numbers}\label{s:elementary}
In this section, we develop a foundation for producing large collections of rational relation numbers in the sequel. 

Let us define
\begin{align*}
R_\bQ&:= \{q\in\bQ\mid  q\text{ is a relation number}\};\\
R_\bQ^{(\ell)}&:=\{q\in\bQ\mid q\text{ is an }\ell\text{--step relation number}\};\\
A_s^{(\ell)}&:=\{r\in\bZ\setminus\{0\}\mid s/r\text{ is an }\ell\text{--step relation number}\}.
\end{align*}
We also let $A_s:=\bigcup_{\ell\ge0} A_s^{(\ell)}$.
\emph{Throughout this section, we fix an integer $s>1$.}
\subsection{On $1$--step relation numbers}

\begin{lem}\label{l:q-basic1}
The following hold.
\be
\item\label{p:1n}
For positive integers $\ell$ and $n$, 
if $q\in R_\bQ^{(\ell)}$, then $\pm q/n\in R_\bQ^{(\ell)}$.
\item\label{p:rst} 
For all nonzero integers $r,s,t$, 
we have $(r+t)/(rst)\in R_\bQ^{(1)}$.
\item\label{p:11n}
For each $n\in\bZ\setminus\{0\}$,
we have that
\[1/n,\ 2/n, 1-1/n\in R_\bQ^{(1)}.\]
\ee
\end{lem}
\bp 
Part (\ref{p:1n}) is immediate from $b_q=(b_{q/n})^n$.
For part (\ref{p:rst}), we let $q=(r+t)/(rst)$ and compute
\[
\xymatrix{
(1,0)\ar@{->>}[r]^r&
\left(1,\frac{r+t}{st}\right)\ar[r]^{-s}&
\left(-\frac{r}t,\frac{r+t}{st}\right)\ar@{->>}[r]^{t}&
\left(-\frac{r}t,0\right).
}\]

Let us prove part (\ref{p:11n}).
Combining Example~\ref{exmp:int} with part (\ref{p:1n}) we see that $1/n$ and $2/n$ are $1$--step relation numbers.  
By substituting $(r,s,t)=(n,1,-1)$, 
we see from part (\ref{p:rst}) that 
\[
1-1/n = -(r+t)/(rst)\]
is  a $1$--step relation number. 
\ep

\subsection{On $2$--step relation numbers}
The notation $x\divides y\pm z$ means $x$ is a divisor of either $y+z$ or $y-z$.
It will be convenient for us to use the notation
 \[
 (x_1,x_2,\ldots,x_k\ ; y ) = \bigcup_{1\le i\le k}(x_i+y\bZ).\]
 For instance, we have $(5\ ; 12 )=5+12\bZ$,
 and  $(\pm5\ ; 12 )=(5+12\bZ)\cup (-5+12\bZ)$.

 The following tool is crucial for this paper. 
\begin{lem}\label{l:quad0}
Suppose there exist nonzero integers $w,m,y$ such that
\[ y\divides m,\quad\text{and}\quad w\divides smy\pm 1.\]
Then for all $r\in (w\ ; sm )\setminus\{0,\pm1,w\}$ we have 
$s/r\in R_\bQ^{(2)}$.
\end{lem}
In particular, it follows that $s/r\in(-4,4)$ for such an $r$.

\bp[Proof of Lemma~\ref{l:quad0}]
We will assume that $w\divides smy-1$, as the other case follows similarly.
For some $u\ne0$ we have
\[1=wu+smy.\]
Let us write $r = w+smt$ for some $t\ne0$, and put $v:=m(y-ut)$. Then 
\[1=(w+smt)u+sm(y-ut)=ru+sv.\]
Since $|r|>1$ and $u\ne0$, we see that $v\ne0$.

After setting $q:=s/r$, we have an orbit of $\form{a,b_q}$ as follows.
\[
\xymatrix{
(1,0)\ar@{>>}[r]^{ry}&(1,sy)\ar[r]^{-v/y}&(ru,sy)\ar@{>>}[r]^>>>>>{-t}&(ru,s(y-ut))\\
{}\ar[r]^>>>>>{m}&(1,s(y-ut))\ar@{>>}[r]^>>>>>{-r(y-ut)}& (1,0).
}
\]
From $rvt\ne0$, it follows that $s/r\in R_s^{(2)}$.\ep


\bd\label{d:sgood}
Let $s,w,m$ be nonzero integers such that $s>1$,
and let \[ D:=\gcd(w,sm),d:=\gcd(w,s).\]
We say the set \[(w\ ; sm )\sse\bZ\]
is an \emph{$s$--good residue class} if
there is an integer $y$  satisfying the following two conditions:
\begin{itemize}
\item $yD\divides md$;
\item $w\divides smy\pm D$.
\end{itemize}
In this case, $w$ is called a \emph{good representative} of $(w\ ; sm )$.
\ed

\begin{exmp}\label{ex:sgood}
\be
\item
The residue class $(0\ ; s )=(s\ ; s )$ is $s$--good.
Indeed, if we set $w=s$ and $m=y=1$, then
\[
w\divides sm-\gcd(w,sm)=0.\]
Moreover, $(0\ ;  sn)=n(0\ ; s )$ is also $s$--good for $n\ne0$.
\item\label{p:splus1}
If $w$ is a divisor of $s\pm 1$, then $(w\ ;s)$ is $s$--good. 
In particular, $\pm(1\ ;s)$ is $s$--good.
\item\label{ex:smy}
More generally, if $w,m,y$ satisfy the hypothesis of Lemma~\ref{l:quad0},
then  $(w\ ; sm )$ is $s$--good. In this case, we have that $\gcd(w,sm)=\gcd(w,s)=1$.
\item
Let $s=25$. If we set $w=9$ and $m=y=2$, then we have \[w=9\divides 99=smy-\gcd(w,sm).\] Hence, $(9\ ;  50)$ is $25$--good.
\ee
\end{exmp}

Recall we have fixed $s>1$ in this section.
We see that all but at most four integers in an $s$--good residue class belong to $A_s^{(2)}$, 
 which generalizes Lemma~\ref{l:quad0}.

\begin{lem}\label{l:quad1}
If $(w\ ; sm )$ is $s$--good with a good representative $w$,
then we have that 
\[
(w\ ; sm )\setminus\{0,w,\pm\gcd(w,sm)\}\sse A_s^{(2)}.\]
\end{lem}


\bp[Proof of Lemma~\ref{l:quad1}]
Let $D$ and $d$ be as in Definition~\ref{d:sgood}.
Set $w'=w/D,s'=s/d$ and $m'=md/D$.
Suppose we have an integer $t$ such that
\[r:=w+smt\not\in\{0,w,\pm D\}.\]
Put $r' :=r/D =w'+ s'm' t$.
By the $s$--good hypothesis, some $y\in \bZ$ satisfies
\[
y\divides m',\text{ and }
w' \divides s'm' y\pm1.\]
Moreover, $r' \not\in\{0,\pm1,w'\}$.
Lemma~\ref{l:quad0} implies that $s'/r'\in R_\bQ^{(2)}$.
It follows that 
\[
\frac{s}{w+smt} = \frac{s'}{r'} \cdot \frac{1}{D/d}\in R_\bQ^{(2)}.\qedhere\]
\ep


Let us note one further consequence of Lemma~\ref{l:quad0}
\begin{lem}\label{l:quad2}
Suppose nonzero integers $w,m,y$ satisfy
\[ y\divides m,\quad\text{and}\quad w\divides smy\pm \gcd(w,s).\]
Then we have that $(w\ ; sm )$ is $s$--good and that
\[
(w\ ; sm )\setminus\{0,w,\pm\gcd(w,s)\}\sse A_s^{(2)}.\]
\end{lem}
\bp
As in Lemma~\ref{l:quad1}, we let $D=\gcd(w,sm)$ and $d=\gcd(w,s)$.
From $D\divides w$ and from the hypothesis, we have $D\divides d$. Indeed, we have \[D\divides w\divides smy\pm d,\]
so that since $D\divides sm$, we have that $d\equiv 0\pmod D$. It follows that $D=d$
and that $(w\ ; sm )$ is $s$--good.
\ep

We note that $z(\{x_i\}\ ; y )=\bigcup_i x_i z+yz\bZ$.
We also record the following.
\begin{lem}\label{l:snm}
If $C$ is an $s$--good residue class, then so is $nC$ for all $n\in\bZ\setminus\{0\}$.
\end{lem}
\bp
Let  $C=(w\ ; sm )$ with a good representative $w$.  Then
$nC=(nw\ ;  snm)$ is also $s$--good; this follows from $\gcd(nw,snm)=|n|\cdot \gcd(w,sm)$.\ep

\begin{prop}
Suppose that for each $n\in\bN$ we can find a collection of  $f(n)$--many $s$--good residue classes whose union contains $\{1,2,\ldots,n\}$. Then we have that 
\[ \underline d\left(A_s^{(2)}\right) \ge 1 - 2\limsup_n f(n)/n.\]
\end{prop}
\bp
By Lemma~\ref{l:quad1}, all positive integers in each $s$--good residue class are in $A_s^{(2)}\cap\bN$, with at most two exceptions.
Hence, we have that
\[ \#(A_s^{(2)}\cap[1,n])/n \ge (n-2 f(n))/n.\qedhere\]\ep

Theorem~\ref{thm:progression} is an immediate consequence of this corollary.

\begin{cor}\label{cor:progression}
For each finite set $Q\sse\bZ$, there is a nonzero integer $M$
such that 
\[Q+sM\left(\bZ\setminus\{0\}\right)\sse A_s^{(2)}.\]
\end{cor}
\bp
For each $w\in Q$, there exists some $m_w\in\bZ\setminus\{0\}$ such that 
\[
sm_w\equiv \gcd(w,s)\pmod{w}.\]
By Lemma~\ref{l:quad2} we have that $(w\ ; sm_w)$ is $w$--good and that
\[w+ sm_w\left(\bZ\setminus\{0\}\right)\sse
A_s^{(2)}\cup \{0,\pm\gcd(w,s)\}.\]

Note that for each $w\in Q$ we have
 \[\gcd(w,s)\in \{d\in\bZ\mid d\text{ divides }s\}.\]
So, for $M_0=\operatorname{lcm}\{m_w\mid w\in Q\}$ we see that
\begin{align*}
Q+sM_0\left(\bZ\setminus\{0\}\right)&\sse
\bigcup\{w+ sM_0\left(\bZ\setminus\{0\}\right)\mid w\in Q\}\\
&\sse
 A_s^{(2)}\cup\{0\}\cup\{d\in\bZ\mid d\text{ divides }s\}.\end{align*}
By setting $M$ to be a sufficiently large multiple of $M_0$, we obtain the desired conclusion.
\ep

\begin{cor}\label{cor:small-denom}
For an integer $w$ in $[-4,4]\cup\{\pm 6\}$, the following hold.
\be
\item\label{p:small-denom1}
The residue class $(w;s)$ is $s$--good.
\item\label{p:small-denom2}
If an integer $t$ satisfies $s/(w+st)\in(-4,4)$, then
$s/(w+st)\in R_\bQ$.
\ee
\end{cor}
\bp
(\ref{p:small-denom1})
By Example~\ref{ex:sgood}, we may only look at the case that $w\ne0$.
It suffices to show that $w$ divides $s\pm \gcd(w,s)$.
We may assume $w\not\divides s$ and $w\not\divides s\pm1$, for otherwise the proof is trivial.
Then it only remains to consider the case $|w|\ge4$.

If $|w|=4$, then our assumption implies that $s\equiv2\pmod{4}$. Then we  see that 
\[s - \gcd(w,s) = s-2\equiv 0\pmod{w}.\]

Suppose $|w|=6$. Our assumption implies that $s\equiv\pm2$ or $s\equiv3$ modulo 6.  Then  $\gcd(w,s)=2$ or $\gcd(w,s)=3$, and we obtain the desired conclusion. 

(\ref{p:small-denom2})
We may assume $w\ne0$. 
Then the above proof implies that $w$ is a good representative of $(w\ ;s)$.
By Lemma~\ref{l:quad1}, we have that either
\[
s/(w+st)\in R_\bQ^{(2)},\]
or
\[w+st\in \{w,\pm\gcd(w,s)\}\sse[-6,6].\]

It is a simple computational verification that for all
nonzero integer $u\in[-6,6]$
and for all integer $s\in(-4|u|,4|u|)$ the number $s/u$ is a relation number; see Proposition~\ref{p:compute} in Appendix~\ref{app:certificate}. This completes the proof that
$s/(w+st)\in R_\bQ$.
\ep
\begin{exmp}
The above corollary implies that $s/(4+st)$ is a $2$--step relation number for all $t\in\bZ$ satisfying $4+st\ne0$ and  $-4<s/(4+st)<4$.
\end{exmp}

The following extends Lemma~\ref{l:q-basic1} (\ref{p:11n}).
\begin{cor}\label{c:basic}
For each nonzero integer $n$,  we have the following:
\[3/n,\ 
1-2/n,\ 
1-3/n,\ 
1-4/n,\ 
1-6/n,\ 
2-1/n
\in R_\bQ\cup\bZ.\]
\end{cor}
\bp
Let $n\in\bZ\setminus\{0\}$ be arbitrary. We may assume $|n|>6$, for otherwise the proof is trivial from direct computations; see also Proposition~\ref{p:compute}.
Since $3$ is a relation number, so is $3/n$. 

In the case when $|w|\le 4$ or $|w|=6$, we see from Corollary~\ref{cor:small-denom} that $1-w/n = (n-w)/ (w+(n-w))$ is a relation number.

Let $w=1-n$ and $s=2n-1$. Since $w\divides s-1$, Lemma~\ref{l:quad0} implies that
\[2-1/n =s/(w+s)\in R_\bQ.\qedhere\]
\ep

\section{Fixed numerators}\label{s:24}
In this section, we establish the Main Conjecture for rational numbers with numerators less than 28 and that are not 24.
\begin{thm}\label{t:sbound}
Let $r,s$ be nonzero integers such that $|s|\le 27$ and $|s|\ne 24$.
If $q=s/r\in(-4,4)$, then $q$ is a relation number.
\end{thm}

We prove Theorem~\ref{t:sbound} for the rest of this section by establishing several claims. 
\emph{We adopt the convention that variables are always integer--valued unless specified otherwise.}

\begin{lem}\label{l:s11}
For each integer $s\in [1,11]\cup\{14,15\}$, we have that
\[
 \{w\in\bZ\mid \gcd(w,s)=1\}=
\bigcup \{(w\ ; s )\mid w\text{ divides }s\pm 1\}.\]
\end{lem}
\bp
If $s=7$ then we see that
\[\{w\in\bZ\mid \gcd(w,7)=1\}=(\pm1,\pm2,\pm3\ ; 7 )=\bigcup\{(w\ ; s )\mid w\text{ divides }6\}.\]
For another example, if $s=11$, then we have 
\[
(\pm1,\pm2,\pm3,\pm4,\pm5\ ; 11 )=\bigcup\{(w\ ; 11 )\in\bZ\mid w\text{ divides }10\text{ or }12\}.\]
The other values of $s$ can be treated similarly, so we omit the details.\ep


\begin{lem}\label{l:nwsm}
Suppose an integer $s$ satisfies $2\le s\le 27$ and $s\ne24$.
\be
\item
Then there exists a finite collection of $s$--good residue classes
\[
\{(w_i\ ; s m_i )\}\]
whose union contains all integers that are relatively prime to $s$;
moreover, we can require that $m_i\divides 60$.
\item
In part (1), we can further require that 
\[\bigcup_i\{w_i,\pm \gcd(w_i,s m_i)\}\sse  A_s\cup[-s/4,s/4].\]
\ee
 \end{lem} 
 
 The requirement that $m_i\divides 60$ in Part (1) of Lemma~\ref{l:nwsm} serves to illustrate the relatively short search that is required
 to find the desired $s$--good residue classes. In order to establish that the set of integers that are relatively prime to $s$ is contained in some union
 of $s$--good residue classes, one may need to exhibit a large number of $s$--good residue classes with moduli which are very big compared to $s$,
 and possibly even unbounded.
 The lemma shows that for small values of $s$ different from $24$, such large moduli are not required.

We note one consequence of Part (2). Suppose $r$ is an integer relatively prime to $s$.
Then  $r$ belongs to $(w_i\ ;  sm_i )$ for some $i$ by Part (1). Lemma~\ref{l:quad1} implies that 
either $s/r$ is a 2--step relation number or 
 \[r\in\{w_i,\pm \gcd(w_i,s m_i)\}.\]
In this latter case, as long as $s/r$ avoids the obvious obstruction that $|s/r|\ge 4$, we will have that $s/r$ is a relation number. 
This point will be crucial in the proof of Theorem~\ref{t:sbound} given at the end of this section. 

\bp[Sketch of the proof of Lemma~\ref{l:nwsm}]
This lemma is a consequence of Proposition~\ref{p:sgood} (1) in Appendix.  For illustration, we will give more hands-on explanation here and leave the computational details to Appendix.

Let us set
\[
X_s:=\left\{
r\in\bZ\mid \gcd(r,s)=1\text{ and }
r\not\equiv w\pmod{s}\text{ for all divisor }w\text{ of }s\pm1\right\}.\]
For part (1), it suffices to find a finite collection of $s$--good residue classes
whose union contains $X_s$;
for, once such a collection is found then we can additionally include
 $(w\ ; s )$ 
for all divisor  $w$ of $s\pm1$. Here, we  are using  Lemma~\ref{l:quad2} in the case $m=1$ and $y=1$.
By Lemma~\ref{l:s11}, we may  assume $s>11$ and $s\not\in\{14,15\}$.

In each case, we will find a list of pairs  $((w\ ;sm),y)$ that satisfy the conditions of Definition~\ref{d:sgood}; we may say $y$ is the ``certificate'' for the $s$--goodness of $(w\ ; sm )$.
We only illustrate the proof for $s=12$ and $s=21$.

\underline{Case $s=12$:} 
Note that $X_s=(\pm5\ ; 12 )=(\pm5,\pm7\ ; 24)$.
Then the following is the desired list of pairs $((w\ ; sm),y)$:
\[ ((\pm5\ ; 24),1),((\pm7\ ; 24),2).\]
This notation is actually an abbreviation of the list
\[ ((5\ ; 24),1),((-5\ ; 24),1),((7\ ; 24),2),((-7\ ; 24),2).\]

\underline{Case $s=21$:}
We have 
$X_s=(\pm8\ ; 21 )$.
We compute as follows.
\begin{align*}
(\pm8\ ; 21 )&=(\pm8,\pm29,\pm13\ ; 63 ),\\
(\pm29\ ; 63 )&=(\pm29,\pm34\ ; 126 )
=
(\pm29\ ; 126 )\cup2(\pm17\ ; 63 ),\\
(\pm13\ ; 63 )&=(\pm13,\pm50\ ; 126 )=(\pm13\ ; 126 )\cup2(\pm25\ ; 63 ),\\
X_s&=
(\pm8\ ; 63 )
\cup
(\pm13,\pm29\ ; 126 )\cup2(\pm17,\pm25\ ; 63 )\\
&\sse
(\pm8\ ; 63 )
\cup
(\pm13,\pm29\ ; 126 )
\cup2(\pm4\ ; 21 )
\end{align*}
Since $(\pm4\ ; 21 )$ is $s$--good, so is $2(\pm4\ ; 21 )$; see Lemma~\ref{l:snm}.
The following is the desired list of pairs:
\[(2(\pm4\ ;s),1),((\pm8\ ;3s),1),((\pm13\ ;6s),3), ((\pm29\ ;6s),3).\]
See Proposition~\ref{p:sgood} for other cases of $s$ and for more details.

For part (2), recall that an $s$--good residue class $(w\ ; sm )$ contains at most three nonzero integers
\[
w,\
\gcd(w,sm),\
-\gcd(w,sm)\] that are possibly not in $A_s^{(2)}$. 
We collect such possible exceptions and individually verify that each one belongs to $A_s$ as long as $|s/r|<4$. This is also done in the proof of Proposition~\ref{p:sgood}.
 \ep

\bp[Proof of Theorem~\ref{t:sbound}]
We may assume that $s>0$.
We have noted after the proof of Lemma~\ref{l:nwsm} that
if $\gcd(r,s)=1$, then $r\in A_s$.

Let us now assume $d:=\gcd(r,s)>1$. Put $r'=r/d$ and \[s'=s/d\le s/2\le 27/2.\]
Since $\gcd(r',s')=1$ and $s'\le 13$, we see from the previous paragraph that $s/r=s'/r'$ is a relation number.
\ep

\section{The case $s=24$}
In this section, we will deduce Theorem~\ref{thm:comput}~(\ref{p:s24}) by proving the following.

\begin{thm}\label{t:s24}
Let $s=24$.
Then there exists a sequence of pairs of integers
\[
\{(a_i,b_i)\}_{i\ge0}\]
such that for each $i\ge0$ and for $M_i=1680\cdot 3^i$,
every integer $x$ satisfies at least one of the following:
\be[(A)]
\item\label{p:remt} We have  $(x\ ;s M_i)\sse(\pm a_i,\pm b_i\ ;{sM_i})$;
\item\label{p:goodt} We have $(x\ ;{sm})$ is an $s$--good residue class for some $m$ dividing $M_i$.\ee
\end{thm}

\bp[Proof of Theorem~\ref{thm:comput}~(\ref{p:s24}) from Theorem~\ref{t:s24}]
Let $s=24$, and let $i\ge0$. 
Recall from Lemma~\ref{l:quad1} that all but at most four integers in each $s$--good residue class belong to $A_s^{(2)}$. Hence, Theorem~\ref{t:s24} implies that
\[
\overline{d}\left(\bZ\setminus A_s^{(2)}\right)
\le 
\overline{d}\left(
(\pm a_i,\pm b_i\ ;{sM_i})\setminus A_s^{(2)}
\right)\le
{d}
(\pm a_i,\pm b_i\ ;{sM_i})\le 4/(sM_i)
.\]
By sending $i\to\infty$, we see that $\bZ\setminus A_s^{(2)}$ has density zero.
\ep

\begin{rem}\label{rem:brute}
A crucial point for the proof of Theorem~\ref{thm:comput}~(\ref{p:s24}) is the choice of  $M_i=1680\cdot 3^i$ as given in Theorem~\ref{t:s24}.
Through a long sequence of trials, errors, and searches by brute force, the authors discovered that the number of non-$s$--good residue classes modulo $sm$ is constant for the choice $s=24$ and $m=M_i$. 
More precisely, the authors use Mathematica to enumerate the number of non-$s$--good residue classes modulo $s t_1 t_2\cdots t_k$ for various choices of $k$ and $t_i\in[2,7]$.
Then it was finally observed that the choices
\[(t_1,t_2,t_3,\ldots) = (4,3,5,7,4,3,3,\ldots)\]
eventually stabilizes the number of non-$s$--good residue classes.
Hence, a suitable modulus to consider is \[4\cdot 3\cdot 5\cdot 7\cdot 4\cdot3\cdot3\cdots =1680\cdot 3^i.\]
As we see below, the justification of this observation will require some arithmetic analysis on the list of non-$s$--good residue classes for each modulus $sM_i$.
\end{rem}

In the remainder of this section, we prove Theorem~\ref{t:s24}. 
A key observation is that the (possibly non-$s$--good) residue classes
\[(\pm a_i,\pm b_i\ ; sM_i)\]
can be expressed by some period--eight sequence $\{z_i\}$.

To be more precise, let us define  integer sequences $\{z_i\}$ and $\{\delta_i\}$ determined by the following conditions.
\begin{itemize}
\item $z_0=1$;
\item $z_i\equiv\delta_i\pmod{3}$ and  $\delta_i\in\{-1,0,1\}$  for each $i\ge0$;
\item $z_{i+1}=(z_i+32\delta_i)/3$.
\end{itemize}

\begin{lem}\label{l:zi}
For each $i\ge0$ and for each $\delta\in\{-1,0,1\}\setminus\{\delta_i\}$,
we have that
\[
z_i+32\delta\sse(\pm1, \pm5\ ;{12}).\]
\end{lem}
\bp
By the nature of the given recursion, the sequences $\{z_i\}$ and $\{\delta_i\}$ must be periodic. 
So, one can verify the lemma by brute force.
Actually, those sequences have period eight; see Table~\ref{t:ci}.
\ep

\begin{table}[hbt!]
\begin{center}
\begin{tabular}{ |p{.3cm}|p{.4cm}|p{3.1cm} |p{1.5cm}||p{.3cm}|p{.4cm}|p{3.1cm}|p{1.5cm}|}
\hline
$i$ & $\delta_i$ & $ z_i,z_i+32,z_i-32$ &   $ z_i+32\delta_i$ & $i$ & $\delta_i$ &  $ z_i,z_i+32,z_i-32$ &  $z_i+32\delta_i$ \\
\hline
\hline
0&1&1,33,-31& 33& 4&0&-15,17,-47&-15\\
\hline
1&-1&11,53,-21& -21& 5&1&-5,27,-37&27\\
\hline
2&-1&-7,25,-39& -39& 6&0&9,41,-23&9\\
\hline
3&-1&-13,19,-45& -45& 7&0&3,35,-29&3\\
\hline
\end{tabular}
\caption{Proof of Lemma~\ref{l:zi}}\label{t:ci}
\end{center}
\end{table}

We can now define the desired sequences $\{a_i,b_i\}_{i\ge0}$ as follows.
\begin{align*}
a_i&=1 + 1260\cdot 3^i z_i,\\
b_i&=1 + 1260\cdot 3^i z_{i+5}. 
\end{align*}

A major computational step of the proof is the following lemma.

\begin{lem}\label{l:ai}
For each $i\ge0$, the following hold.
\be
\item\label{p:ai1}
$a_{i+1}=a_i+\delta_i s  M_i$
and $b_{i+1}=b_i+\delta_{i+5} s  M_i$
\item\label{p:del1} If \[\delta\in\{-1,0,1\}\setminus\{\delta_i\}\]
then 
there exists a divisor $m$ of $M_{i+1}$ such that 
\[
(a_i+\delta s  M_i\ ;{sm})\]
is $s$--good.
\item\label{p:del2} If \[\delta\in\{-1,0,1\}\setminus\{\delta_{i+5}\}\]
then 
there exists a divisor $m$ of $M_{i+1}$ such that 
\[
(b_i+\delta s  M_i\ ;{sm})\]
is $s$--good.
\ee
\end{lem}
\bp
(\ref{p:ai1})
Note that
\[1260=35\cdot 36,\quad 1680s=35\cdot 36\cdot 32.\]
We see from the definitions of $\{z_i\}$ and $\{\delta_i\}$ preceding Lemma~\ref{l:zi} that
\begin{align*}
&a_i+\delta_i s M_i-a_{i+1}
= 1260\cdot 3^i z_i + 1680s\cdot 3^i \delta_i -1260\cdot3^{i+1}  z_{i+1}\\
&=35\cdot 36\cdot 3^i ( z_i + 32\delta_i-3z_{i+1})=0.\\
&b_i+\delta_{i+5} s M_i-b_{i+1}
= 1260\cdot 3^i z_{i+5} + 1680s\cdot 3^i \delta_{i+5} -1260\cdot3^{i+1}  z_{i+6}=0.
\end{align*}

(\ref{p:del1}) We saw in Lemma~\ref{l:zi} that 
\[ z_i+32\delta= 12p+ c\]
 for some $p\in\bZ$ and $c\in\{\pm1,\pm5\}$.
Put $m= 18c\cdot 3^i$, so that $m\divides M_{i+1}$. Then 
\begin{align*}
&(a_i+\delta s  M_i-1+2 m)/(sm)
=3^i(1260  z_i+1680s \delta   + 36c)/(3^i\cdot36 \cdot  12c)\\
&=( 35(  z_i+32 \delta ) + c)/(12c)
= (35\cdot(12p+c)+c)/(12c) =(35/c)p+3\in\bZ.
\end{align*}
So, we have that
\[
a_i+\delta s  M_i\equiv 1-2m\pmod{sm}.\]
Note that
\[
1-2m\divides 4m^2-1=s \cdot m\cdot (m/6)-1.\]
By setting $y=m/6$ in Definition~\ref{d:sgood} (or, Example~\ref{ex:sgood} (\ref{ex:smy})), we see that
\[
(a_i+\delta s  M_i\ ;{sm})=(1-2m\ ;{sm})\]
is an $s$--good residue class.

(\ref{p:del2}) The proof is essentially the same, after replacing $(a_i,\delta_i)$ by $(b_i,\delta_{i+5})$. 
\ep

\bp[Proof of Theorem~\ref{t:s24}]
We use induction. The base case $i=0$ is a consequence of Proposition~\ref{p:sgood} in Appendix, where a computer--assisted proof is given.  Namely, we may set \[a_0=1261,\quad b_0=-6299.\]

Let us now assume the conclusion for some $i\ge0$. 
To obtain a contradiction, we also assume that neither of the alternatives (\ref{p:remt}) or (\ref{p:goodt}) holds for the index $i$ and for some fixed positive integer $x$. 

In the case when 
\[x\not\in (\pm a_i,\pm b_i\ ;{sM_i}),\]
we see from the inductive hypothesis that $(x\ ;{sm})$ is $s$--good for some $m\divides M_i$.
Since
$M_i\divides M_{i+1}$, the alternative (\ref{p:remt}) holds for the index $i+1$, we are done with this case.

We will now consider the case that
\[
x\in (\pm a_i,\pm b_i\ ;{sM_i}).\]
Let us first suppose 
\[x\in  (a_i\ ; sM_i)= (a_i- sM_i, a_i, a_i+ sM_i\ ;{sM_{i+1}}).\]
Then we have $x\in (a_i + s\delta M_i\ ; sM_{i+1})$ for some $\delta\in\{-1,0,1\}$.
If $\delta=\delta_i$, then Lemma~\ref{l:ai} implies that
\[x\in (a_i+\delta_i s M_i\ ;{sM_{i+1}})=(a_{i+1}\ ;{sM_{i+1}}),\]
and that the alternative (\ref{p:remt}) for the index $i+1$ is satisfied.
If $\delta\ne\delta_i$, then 
the same lemma implies that $x$ satisfies the alternative (\ref{p:goodt}) for the index $i+1$. This completes the proof for the case $x\in  (a_i\ ; sM_i)$.

By applying the same argument to the residue classes
\[-(a_i\ ;sM_i), \ (b_i\ ;sM_i),\  -(a_i\ ;sM_i)\]
we obtain the desired conclusion for $i+1$.
\ep


\section{General density estimates}\label{sec:gde}
In this section we establish Theorem~\ref{thm:density}, which we do by an averaging argument. The general strategy is as follows:
suppose $X\sse\bN\times\bN$, with horizontal and vertical sections $H_i=\{y\mid (y,i)\in X\}$ and $V_i=\{z\mid (i,z)\in X\}$ respectively.
One is interested in estimating the density in $\bN$ of the horizontal sections $H_i$ of $X$ from below, but these may be difficult to compute.
However, one may have better methods for computing the vertical sections $V_i$ of $X$. So, one truncates $\bN\times\bN$ to $[1,H]\times [1,V]$ for some suitably chosen large values of $H$ and $V$, and one
adds up the sizes of the vertical sections $V_i$ of $X$ restricted to $i=[1,H]$. Dividing by $V$ gives the average size of a vertical section
of $X$.

In more specific terms, we fix a numerator $s$ and a large multiple $sm$ of $s$, which serves as the truncation $V$ above. One then enumerates residue classes modulo $sm$ which are not contained in $s$--good residue classes (subject to some further constraints to make calculations more tractable), and the number $\ell$ of these serves as the truncation $H$. The number--theoretic lemmata developed earlier allow us to then estimate the density of non--relation numbers of the form $s/r$. We now make this approach precise.

Fix $s\ge 28$.
In what follows, we recursively construct an increasing sequence $\{m_n\}_{n\ge0}$
such that the set
\[
B_n:=\left\{
(w\ ; sm_n)\middle\vert
(w\ ;sm_n)\text{ is not contained in an }s\text{--good residue class}\right\}\]
has a small density. Then, we apply Lemma~\ref{l:quad1} to see that
\[\underline d\left(A_s^{(2)}\right)\ge 1- \# B_n/(sm_n).\]

We begin by setting $m_0=1$. By Corollary~\ref{cor:small-denom}, 
we see that \[(w\ ;s)\not\in B_0\] for $|w|\le4$ or $|w|=6$. In particular,
\[
\# B_0/ (sm_0)\le 1-11/s.\]

Suppose we have constructed $m_n\in\bN$. For brevity, let us write
\[ m:=m_n,\ B_n=\{(w_1\ ; sm), \ldots,(w_\ell\ ; sm)\},\ 
v_i:=\gcd(w_i,sm)\]
for some $\ell>0$.
We may choose $w_i$ in the set $(-sm/2,sm/2]$ such that \[w_i\not\in \{0,\pm1,\pm2,\pm3,\pm4,\pm6\}.\]
We define
\[
Z:=\{(i,x)\mid i\in[1,\ell]\text{ and }x\in[1,sm]\text{ such that }smx\equiv \pm v_i\pmod{w_i}\}.\]
Let $Y_i:=Z\cap (\{i\}\times\bZ)$.
We begin by establishing the following.

\begin{claim}\label{cl:zset}
The following hold.
\be
\item\label{p:zgood}
For each $(i,x)\in Z$ the residue class $(w_i\ ;smx)$ is $s$--good.
\item\label{p:zdiff}
If $(i,x)$ and $(j,x)$ are distinct elements of $Z$, then the residue classes
$(w_i\ ;smx)$ and $(w_{j}\ ;smx)$ are distinct as well.
\item\label{p:zcount} 
For each $i\in[1,\ell]$, the cardinality of $Y_i$ is at least four.
\ee
\end{claim}
\bp[Proof of Claim~\ref{cl:zset}]
(1) If $(i,x)\in Z$, then $\gcd(w_i,smx)=v_i$. Applying Definition~\ref{d:sgood} after replacing $m$ by $mx$, setting $y=1$, and writing
$D=v_i$,
we have $w_i\divides smx\pm v_i$, whence we may conclude that $(w_i\ ;smx)$ is $s$--good.

(2) This is becase $w_i\not\equiv w_j\pmod{sm}$.

(3) Suppose first that $w_i\not\divides 2v_i$. 
Then the modular arithmetic equation
\[ smx\equiv  v_i\pmod{w_i}\]
has a unique solution modulo $w_i/v_i$.
Similarly, $smx\equiv-v_i$ has a unique solution as well.
Since $v_i\not\equiv -v_i\pmod{w_i}$, we have that
\[\#Y_i \ge 2\floor*{\frac{sm}{|w_i/v_i|}}
\ge 2\floor*{\frac{sm}{sm/2}}=4.\]

If $w_i\divides 2v_i$, then $w_i/v_i=\pm2$ or $\pm1$. So, we have
\[\#Y_i \ge \floor*{\frac{sm}{2}}
\ge \frac{28}{2}>4.\qedhere\]
\ep

By  applying Claim~\ref{cl:zset}
and 
averaging $x$ over $[1,\ell]$, we can find some $X\in[1,sm]$ such that
the number of  distinct $s$--good residue classes in the set \[\{(w_i\ ;smX)\}_{i\in[1,\ell]}\]
 is at least 
\[
\#Z/sm\ge \sum_i\# Y_i / sm \ge 4\ell/sm.\]
To make the recursion deterministic, we pick the smallest such $X$. 

We now define $x_n:=X$ and $m_{n+1}=m_nx_n$. 
The set $B_{n+1}$ is contained in the set 
\[\{(w_i+smk\ ;smX)\mid i\in[1,\ell]\text{ and }k\in[0,X)\}.\]
In the set above, at least $4\ell/sm$ residue classes are $s$--good.
It follows that
\[
\frac{\# B_{n+1}}{sm_{n+1}}
\le 
\frac{1}{smX}
\left(
\ell X - \frac{4\ell}{sm}
\right)
=
\frac{\# B_n}{sm_n}
\left(
1 - \frac{4}{sm_{n+1}}
\right)
\]

Summing up, we have that
\[
\bar d\left(\bZ\setminus A_s^{(2)}\right)
\le
\liminf_{n\to\infty} \frac{\# B_n}{sm_n}
\le
\left(1-\frac{11}{s}\right)
\prod_{n=1}^\infty \left(1-\frac{4}{sm_n}\right).
\]

From the inequality $x_{n-1}\le sm_{n-1}$, we have that 
\[
m_n=sm_{n-1} x_{n-1}\le s^2 m_{n-1}^2\le \cdots
\le
s^{2+4+\cdots+2^n}m_0^{2^n}=s^{2^n-2}.\]
Hence, the theorem follows. \qed

As remarked in the introduction, Theorem~\ref{thm:density} does not quite show that $A_s^{(2)}$ has natural density $1$, but the infinite
product does give a significant improvement to the density estimate. As a particular example, we consider the case $s=28$. We have
that
\[\left(1-\frac{11}{28}\right)=\frac{17}{28}\approx 0.6071428571.\] The infinite product converges very quickly, and multiplying it out up to $n=4$
yields
\[\left(1-\frac{11}{28}\right)\left(1-\frac{4}{28}\right)\left(1-\frac{4}{28^3}\right)\left(1-\frac{4}{28^7}\right)\left(1-\frac{4}{28^{15}}\right)\approx
0.5203133366.\]

Similarly, for $s=29$ we obtain the estimates $0.6206896552$ and $0.5349895317$, respectively. For $s=30$, we obtain the estimates
$0.6333333333$ and $0.5488075719$, respectively.

\appendix
\section{Certifying $s/r$ is a relation number}\label{app:certificate}
In this appendix, we give a detailed description of the algorithms used in the paper.
The Mathematica code implementing such algorithms, as well as the relevant outputs of those code, are available for download as an ancillary file  (\texttt{relnum-v2.pdf}) with the arXiv version of this paper~\cite{KK2019:compute} and also on the authors' respective websites.

Conceptually, for each orbit point $(x_i,y_i)$  Algorithm~\ref{a:relnum}  determines the next orbit point $(x_{i+1},y_{i+1})$, so that either $|x_{i+1}|$ or $|y_{i+1}|$ is minimized (depending on the parity of $i$) over all possible choices. This algorithm is inspired by~\cite{LyndonUllman69,TanTan96}.

To be more precise, let $x,y\in\bZ\setminus\{0\}$. 
Setting $t=x-y\floor{x/y}$,
we define a \emph{shifted remainder} of $x$ by $y$ as
\[
\textsc{SR}(x,y):=
\begin{cases}
t,&\text{ if }t=x+y,\\
t,&\text{ if }t\ne x\text{ and }|t|\le |y|/2,\\
t - y,&\text{ otherwise}\end{cases}\]
Note that 
\[|\textsc{SR}(x,y)|=\min\{ |x+yk| \co k\ne0\}.\]
We also let
\[
\operatorname{sign}(x):=
\begin{cases}
1,&\text{ if }x\ge0,\\
-1,&\text{ otherwise}.\end{cases}\]

The function $\textsc{RelNum}(s,r,M)$ given in 
Algorithm~\ref{a:relnum} can determine (when it succeeds) that a given number $s/r$ is a relation number under $M$ iterations. 
This algorithm begins with the moves
\begin{align*}
(1,0)&=(r,0)\stackrel{1}\tto (r,s)\to (\textsc{SR}(r,s),s)=  (
r\textsc{SR}(r,s),rs)\\
&\stackrel{-1}\tto 
(r\textsc{SR}(r,s),s(r-\textsc{SR}(r,s) )=:(x_0,sy_0).
\end{align*}
Using the variables $d=\gcd(x'',y')$ and $\sigma=\operatorname{sign}(x'')$, we then define
\begin{align*}
&(x_i,sy_i)\to (\textsc{SR}(x_i, sy_i),sy_i)=(x',sy_i)
= (rx',rsy_i)
\tto(rx',s\cdot \textsc{SR}(ry_i, x'))\\
&=(x'',sy')=(x''\sigma/d,sy'\sigma/d)=(x_{i+1},y_{i+1}).\end{align*}

The function $\textsc{RelNum}(s,r,M)$ returns \texttt{True}
 if the orbits $\{(x_i,y_i)\}$ becomes periodic (up to changing the sign of $y_i$), or if $x_iy_i(x_i-1)=0$ for some $i\le M$.
In this case, we see that $s/r$ is a relation number; see Proposition~\ref{p:CF}.
Otherwise, the algorithm returns \texttt{False}, and is inconclusive.

\begin{algorithm}
\caption{Certifying $q=s/r\in R_\bQ$ by shifted remainders}\label{a:relnum}
\begin{algorithmic}[1]
\Function{RelNum}{$s,r,M$}
	\If {$|s/r|\ge4$ or $\gcd(s,r)\ne1$ or $r\in\bZ$} 
		\State{\textbf{Print}(``known cases") and \Return{\texttt{Null}}}
	\EndIf
	\State $ i\gets 0,\quad x_0\gets r,\quad y_0\gets (r-\textsc{SR}(r, s))/s,\quad \text{flag} \gets \texttt{False}$
	\If {$y_0=0$} flag$\gets{\texttt{True}}$\EndIf
	\While{$i<M$ and flag$=$\texttt{False}}
		\State{$x\gets x_i$ and $y\gets y_i$}
		\State{$x'\gets \textsc{SR}(x,sy)$}
		\State{$y'\gets \textsc{SR}(ry,x')$ and $x''\gets rx'$}
		\State{$d\gets \gcd(x'',y')$ and $\sigma=\operatorname{sign}(x'')$}
		\If{d$\ne0$}
			\State{$x\gets x''\sigma/d$ and $y\gets y'\sigma/d$}
		\Else
			\State{$x\gets x''\sigma$ and $y\gets y'\sigma$}			
		\EndIf		
		\If{$xy(x-1)=0$ or $(x,y)= (x_j,\pm y_j)$ for $\exists j<i$}
				\State{flag$\gets$\texttt{True}}
		\EndIf
		\State{$i\gets i+1$}
		\State{$x_i\gets x$ and $y_i\gets y$}
	\EndWhile
	\State{\Return{flag}}
\EndFunction
\end{algorithmic}
\end{algorithm}

Let us now consider a (typically slower) variation of Algorithm~\ref{a:relnum}. 
Again, for a given $(x_i,y_i)$  this algorithm tries to find a sequence
\[
(x_i,y_i)\tto(x_{i+1},y_{i+1})\to(x_{i+2},y_{i+2})\]
so that $|x_{i+2}|$ is minimized among possible choices. 

To be precise, for nonzero integers $a,b,c$ satisfying $a,c>0$
and $a\not\divides b$, we write
\[(u,v)=\textsc{Min}_3(a,b,c),\]
where $u$ and $v$ are nonzero integers minimizing the value
\[|(au+b)v+c|.\]
We consider an arbitrary choice if such a pair $(u,v)$ is not unique.
Then Algorithm~\ref{a:relnummin} attempts to find an orbit coming from the moves
\[
(x,sy)=(rx,rsy)\tto(rx,s(ry+ux))\to(rx+s(ry+ux)v,s(ry+ux))=(x',sy'),\]
while minimizing the value of $|x'|= | (sx\cdot u+sry)\cdot v+rx|$ in each step
by setting
\[
(u,v) = \textsc{Min}_3(sx,sry,rx).\]
So, $\textsc{RelNumMin}(s,r,M)$ functions exactly as $\textsc{RelNum}(s,r,M)$, except that it uses Algorithm~\ref{a:relnummin}.

\begin{algorithm}
\caption{Certifying $q=s/r\in R_\bQ$ by minimizing coordinates}\label{a:relnummin}
\begin{algorithmic}[1]
\Function{RelNumMin}{$s,r,M$}
	\If {$|s/r|\ge4$ or $\gcd(s,r)\ne1$ or $r\in\bZ$} 
		\State{\textbf{Print}(``known cases") and \Return{\texttt{Null}}}
	\EndIf
	\State $ i\gets 0,\quad x\gets \textsc{SR}(r,s),\quad y\gets (r- \textsc{SR}(r,s))/s,\quad \text{flag} \gets \texttt{False}$
	\If {$y_0=0$} flag$\gets{\texttt{True}}$\EndIf
	\State{$x_0\gets \operatorname{sign}(x)x$ and $y_0\gets  \operatorname{sign}(x)y$}
	\While{$i<M$ and flag$=$\texttt{False}}
		\State{$x\gets x_i$ and $y\gets y_i$}
		\State{$(u,v)\gets \textsc{Min}_3(sx,sry,rx)$}
		\State{$x'\gets rx+s(r y+ux)v$ and $y'\gets ry+ux$}
		\State{$d\gets \gcd(x',y')$ and $\sigma=\operatorname{sign}(x')$}		
		\If{d$\ne0$}
			\State{$x\gets x'\sigma/d$ and $y\gets y'\sigma/d$}
		\Else
			\State{$x\gets x'\sigma$ and $y\gets y'\sigma$}
		\EndIf
		\If{$xy(x-1)=0$ or $(x,y)= (x_j,\pm y_j)$ for $\exists j<i$}
				\State{flag$\gets$\texttt{True}}
		\EndIf
		\State{$i\gets i+1$}
		\State{$x_i\gets x$ and $y_i\gets y$}
	\EndWhile
	\State{\Return{flag}}
\EndFunction
\end{algorithmic}
\end{algorithm}

The following conjecture would imply the Main Conjecture.
\begin{con}
For all $s,r\in\bN$ satisfying $s/r<4$,
there exists $M>0$ such that $\textsc{RelNum}(s,r,M)=\texttt{True}$
or
$\textsc{RelNumMin}(s,r,M)=\texttt{True}$.
\end{con}

Using these algorithms, we prove the following.
\begin{prop}\label{p:compute}
Let $s$ and $r$ be positive integers such that $s/r<4$.
\be
\item If $r\le 8$, then $s/r$ is a relation number.
\item If $s\le 30$ and $s/r\ge 1/10$, then $s/r$ is a relation number.
\ee
\end{prop}
\bp
By induction, it suffices to consider the case when $\gcd(s,r)=1$.

For part (1), we use the Mathematica to compute the value $\textsc{RelNum}(s,r,5000)$ for each $2\le r\le 8$ and $2\le s\le 4r-1$.  The result shows that all rational numbers $s/r$ in this range are relation numbers, and that this can be verified under 5000 iterations with Algorithm~\ref{a:relnum}. We also remark that 
\[ 35/9,\ 39/10\]
 are inconclusive under 5000 iterations.

For part (2), we again apply the function $\textsc{RelNum}(s,r,5000)$ for each $2\le s\le 30$ and $2\le r\le 10s$. The output says that all rational numbers $s/r$ in this range are relation numbers possibly except for 
\[
28/17,29/17.\]

For the above two rational numbers, we then apply Algorithm~\ref{a:relnummin}. 
The output of the second algorithm then tells us that these two numbers are indeed relation numbers. For instance, when $s/r = 29/17$ this second algorithm finds a sequence
\begin{align*} &(17,29) = (-17,-29) \to(12,-29)=(12\cdot 17,-29\cdot 17) \\
&\tto(12\cdot 17,29(-17+12\cdot2))=(204,203)\to(1,203).\end{align*}
So, we are done.
\ep

\section{Certifying $\bZ$ is a finite union of $s$--good residue classes}
\begin{prop}\label{p:sgood}
Let $s$ be a positive integer in $[2,27]$.
\be
\item\label{p:sn24b}
If $s\ne 24$, then there exists a finite collection of $s$--good residue classes
\[
\{(w_i\ ; sm_i)\}_{1\le i\le k}\]
whose union is $\bZ$, such that \begin{equation}\label{eq:smi}\tag{**}
\bigcup_{1\le i\le k}\{w_i, \pm\gcd(w_i,sm_i)\}\sse A_s\cup[-s/4,s/4].\end{equation}
Moreover, we can require that $m_i\divides 60$.
\item\label{p:s24b}
If $s=24$, then every integer $x$ satisfies at least one of the following.
\be[(A)]
\item\label{p:rem} we have that  $(x\ ;1680s)\sse(\pm 1261,\pm 6299\ ;{1680s})$;
\item\label{p:good} we have that $(x\ ;{sm})$ is an $s$--good residue class for some $m$ dividing $1680$.\ee
\ee
\end{prop}
\bp
(1)
By induction, it suffices to find a finite collection $\FF_s=\{(w_i\ ; sm_i)\}_i$ of $s$--good residue classes containing
\[Y_s=\{x\in\bN\mid \gcd(x,s)=1\}\]
such that (\ref{eq:smi}) holds.

If $s\le11$, then we simply choose the collection
\[\FF_s=\{(x\ ; s)\mid \gcd(x,s)=1\text{ and }1\le x<s\}.\]
From Lemma~\ref{l:s11}, each residue class in the above collection is $s$--good.
Moreover, whenever $1\le x<s$ we have that $x\in A_s$ by Proposition~\ref{p:compute}.
This completes the proof for $s\le 11$.

Let $s\ge12$. Let us list a specific sequence $t_{12},t_{13},\ldots$ as follows.
\[
t_{12}=2,2,2,2,2,3,6,2,4,6,12,12,1680,6,60,t_{27}=60.\]
In particular, $t_{24}=1680$; see Remark~\ref{rem:brute} regarding the choice of $t_{24}$.

Except for the case $s=24$, this sequence $\{t_s\}$ is found by brute force in the range $t_s\in [2,60]$ until the set $\bZ$ is completely covered by $s$--good residue classes.
More precisely, the number $t_s$ satisfies the following claim.

\begin{claim}
Let $s\in[12,27]$ and $s\ne24$. Then for each $w\in [1,st_s)$ satisfying $\gcd(w,s)=1$,
there exist integers $w',m,y$ such that the following hold:
\be[(i)]
\item $y\divides m$ and $m\divides t_s$;
\item $w'\equiv \pm w\pmod{sm}$;
\item $smy\equiv\pm\gcd(w',sm)\pmod{w'}$;
\item $w',\gcd(w',sm)\in A_s\cup[-s/4,s/4]$.
\ee
\end{claim}

The claim again can be proved by a brute force search, as illustrated in the ancillary file.
This search is successful in the finite range $|w'|\le |w|$ and $y\divides m$ and $m\divides t_s$.

Once the claim is proved, note that 
Parts (i) through (iii), along with Lemma~\ref{l:quad1}, imply each element $r$ in the residue class $(w\ ; s t_s)$ belongs to $A_s^{(2)}$ with possible exceptions of
\[
r\in \{0,w',\pm\gcd(w',sm)\}.\]
In these exceptional cases, Parts (iv) implies that $r\in A_s$ unless $|s/r|>4$.
In particular,  Part (\ref{p:sn24b}) is proved.

For Part (\ref{p:s24b}), we again run the same algorithm for $s=24$ as in Part (1). We then observe that Parts (i) through (iv) of the above claim holds as long as 
\[w\not\in \{1261,6299,34021,39059\}.\]
This implies Part (\ref{p:s24b}).\ep

\section*{Acknowledgements}
The authors thank an anonymous referee for helpful comments and corrections.
The authors thank V. Shpilrain for pointing out the Main Conjecture to them.
The authors also thank T. Tsuboi for suggesting a connection to generalized continued fractions. 
The first author is supported by Samsung Science and Technology Foundation (SSTF-BA1301-51) and by a KIAS Individual Grant (MG073601) at Korea Institute for Advanced Study.
The second author is partially supported by an Alfred P. Sloan Foundation Research Fellowship and NSF Grant DMS-1711488.

\def\cprime{$'$} \def\soft#1{\leavevmode\setbox0=\hbox{h}\dimen7=\ht0\advance
  \dimen7 by-1ex\relax\if t#1\relax\rlap{\raise.6\dimen7
  \hbox{\kern.3ex\char'47}}#1\relax\else\if T#1\relax
  \rlap{\raise.5\dimen7\hbox{\kern1.3ex\char'47}}#1\relax \else\if
  d#1\relax\rlap{\raise.5\dimen7\hbox{\kern.9ex \char'47}}#1\relax\else\if
  D#1\relax\rlap{\raise.5\dimen7 \hbox{\kern1.4ex\char'47}}#1\relax\else\if
  l#1\relax \rlap{\raise.5\dimen7\hbox{\kern.4ex\char'47}}#1\relax \else\if
  L#1\relax\rlap{\raise.5\dimen7\hbox{\kern.7ex
  \char'47}}#1\relax\else\message{accent \string\soft \space #1 not
  defined!}#1\relax\fi\fi\fi\fi\fi\fi}
\providecommand{\bysame}{\leavevmode\hbox to3em{\hrulefill}\thinspace}
\providecommand{\MR}{\relax\ifhmode\unskip\space\fi MR }
\providecommand{\MRhref}[2]{%
  \href{http://www.ams.org/mathscinet-getitem?mr=#1}{#2}
}
\providecommand{\href}[2]{#2}

\end{document}